\documentclass[12pt,leqno,draft]{article}

\begin{document}

\title{Metric spaces: \\ The definition, and some examples}

\author{Stephen Semmes \\
        Rice University}

\date{}

\maketitle


\renewcommand{\thefootnote}{}   

\footnotetext{These informal notes have been prepared in connection
with a lecture at the high school mathematics tournament held at Rice
University on February 20, 2010.}

        A \emph{metric space} is a set $M$ together with a real-valued
function $d(x, y)$ defined for $x, y \in M$ that satisfies the
following three conditions.  First, $d(x, y) \ge 0$ for every $x, y
\in M$, and $d(x, y) = 0$ if and only if $x = y$.  Second,
\begin{equation}
        d(x, y) = d(y, x)
\end{equation}
for every $x, y \in M$.  Third,
\begin{equation}
        d(x, z) \le d(x, y) + d(y, z)
\end{equation}
for every $x, y, z \in M$, which is known as the \emph{triangle inequality}.
This function $d(x, y)$ is called the \emph{metric} on $M$, and represents
a measurement of distance between elements of $M$.

        Remember that the \emph{absolute value} $|x|$ of a real number
$x$ is equal to $x$ when $x \ge 0$, and is equal to $-x$ when $x \le
0$.  As usual,
\begin{equation}
        |x \, y| = |x| \, |y|
\end{equation}
and
\begin{equation}
        |x + y| \le |x| + |y|
\end{equation}
for every $x, y \in {\bf R}$.  The standard metric on the real line
is defined by
\begin{equation}
        d(x, y) = |x - y|,
\end{equation}
which is easily seen to satisfy the requirements of metric mentioned
before.

        The standard Euclidean metric on ${\bf R}^n$ is defined by
\begin{equation}
        d_2(x, y) = \Big(\sum_{j = 1}^n (x_j - y_j)^2\Big)^{1/2},
\end{equation}
where $x = (x_1, \ldots, x_n)$, $y = (y_1, \ldots, y_n)$.  This
clearly satisfies the positivity and symmetry conditions for a metric,
but the triangle inequality is less obvious.  The latter is a
well-known theorem in classical Euclidean geometry, and can also be
shown using the Cauchy--Schwarz inequality.  Alternatively, it is easy
to check that
\begin{equation}
        d_1(x, y) = \sum_{j = 1}^n |x_j - y_j|
\end{equation}
and
\begin{equation}
        d_\infty(x, y) = \max_{1 \le j \le n} |x_j - y_j|
\end{equation}
also define metrics on ${\bf R}^n$, using the triangle inequality for
the standard metric on the real line.  Although these three metrics on
${\bf R}^n$ are different geometrically, they are equivalent
topologically, which basically means that they are the same in terms
of continuity and convergence.

        If $(M, d(x, y))$ is any metric space, then the (open) ball
$B(x, r)$ with center $x \in M$ and radius $r > 0$ is defined by
\begin{equation}
        B(x, r) = \{x \in M : d(x, y) < r\}.
\end{equation}
In the case of ${\bf R}^n$ equipped with the standard Euclidean
metric, $B(x, r)$ is an ordinary round ball.  However, if ${\bf R}^n$
is equipped with the metric $d_1(x, y)$, then $B(x, r)$ is
diamond-shaped.  If ${\bf R}^n$ is equipped with the metric
$d_\infty(x, y)$, then $B(x, r)$ is a cube with sides parallel to the
coordinate axes.

        If $(M, d(x, y))$ is any metric space and $0 < \alpha \le 1$,
then one can show that $d(x, y)^\alpha$ also defines a metric on $M$.
The main point is to check that the triangle inequality still holds,
by verifying that
\begin{equation}
\label{(a + b)^alpha le a^alpha + b^alpha}
        (a + b)^\alpha \le a^\alpha + b^\alpha
\end{equation}
for every pair $a$, $b$ of nonnegative real numbers.  Equivalently,
\begin{equation}
        a + b \le (a^\alpha + b^\alpha)^{1/\alpha}
\end{equation}
for every $a, b \ge 0$.  If $1/\alpha$ is an integer, then this
follows by expanding the right side into a sum.  Otherwise, (\ref{(a +
b)^alpha le a^alpha + b^alpha}) can be seen as a nice exercise in
calculus.  One can also use
\begin{equation}
        \max (a, b) \le (a^\alpha + b^\alpha)^{1/\alpha}
\end{equation}
to get that
\begin{equation}
        a + b \le (a^\alpha + b^\alpha) \, \max(a, b)^{1 - \alpha}
               \le (a^\alpha + b^\alpha)^{1/\alpha}.
\end{equation}
As before, $d(x, y)^\alpha$ is equivalent to $d(x, y)$ topologically,
even if they may be different geometrically.

        If $M$ is any set, then one can define a metric on $M$ by
putting $d(x, y) = 1$ when $x \ne y$, and equal to $0$ when $x = y$.
It is easy to see that this satisfies the requirements of a metric,
which is known as the \emph{discrete metric} on $M$.

        Let $p$ be a prime number, and let $x = (a/b) \, p^j$ be a
rational number, where $a$, $b$, $j$ are integers, $a, b \ne 0$, and
$a$, $b$ are not divisible by $p$.  The \emph{$p$-adic absolute value}
of $x$ is denoted $|x|_p$ and defined to be $p^{-j}$ under these
conditions, and to be $0$ when $x = 0$.  One can check that
\begin{equation}
        |x \, y|_p = |x|_p \, |y|_p
\end{equation}
and
\begin{equation}
        |x + y|_p \le \max(|x|_p, |y|_p) \le |x|_p + |y|_p
\end{equation}
for every pair of rational numbers $x$, $y$.  The \emph{$p$-adic metric}
on the set ${\bf Q}$ of rational numbers is given by
\begin{equation}
        d_p(x, y) = |x - y|_p.
\end{equation}
In particular,
\begin{equation}
        d_p(x, z) \le \max(d_p(x, y), d_p(y, z)) \le d_p(x, y) + d_p(y, z)
\end{equation}
for every $x, y, z \in {\bf Q}$.  This is quite different from the
standard metric on ${\bf Q}$, since $|x|_p \le 1$ for every integer
$x$, and $|x|_p$ may be quite small even when $x$ is a nonzero
integer.  The \emph{$p$-adic numbers} are obtained by completing the
rational numbers with respect to the $p$-adic metric, in the same way
that the real numbers can be obtained by completing the rational
numbers with respect to the standard metric.

        If $(M, d(x, y))$ is any metric space and $E$ is a subset of
$M$, then the restriction of $d(x, y)$ to $x, y \in E$ is a metric on
$E$.  In this way, every subset of a metric space may be considered as
a metric space too.

        Let ${\bf S}^n$ be the unit sphere in ${\bf R}^{n + 1}$ with
respect to the standard Euclidean metric, so that
\begin{equation}
\label{{bf S}^n = ...}
 {\bf S}^n = \bigg\{(x_1, \ldots, x_{n + 1}) \in {\bf R}^{n + 1} :
                         \sum_{j = 1}^{n + 1} x_j^2 = 1\bigg\}.
\end{equation}
If $x, y \in {\bf S}^n$, $x \ne y$, and $x \ne -y$, then there is a
unique $2$-dimensional plane $P(x, y)$ in ${\bf R}^{n + 1}$ that
passes through $x$, $y$, and $0$.  Thus
\begin{equation}
\label{C(x, y) = P(x, y) cap {bf S}^n}
        C(x, y) = P(x, y) \cap {\bf S}^n
\end{equation}
is a circle with radius $1$, and we let $d_{{\bf S}^n}(x, y)$ be the
length of the shorter arc on $C(x, y)$ that connects $x$ to $y$.  If
$x = y$, then we simply put $d_{{\bf S}^n}(x, y) = 0$.  If $x = -y$,
then $x$, $y$ lie on the same line through $0$, and we can use any
plane $P(x, y)$ that contains this line.  In this case, all of the
circular arcs on ${\bf S}^n$ that connect $x$ to $y$ are half-circles,
and we take $d_{{\bf S}^n}(x, y)$ to be their common length, which is
$\pi$.  It is well known that $d_{{\bf S}^n}(x, y)$ satisfies the
triangle inequality, and this will be discussed further later on.
Using this, it is easy to see that $d_{{\bf S}^n}(x, y)$ defines a
metric on ${\bf S}^n$, and that this metric is topologically
equivalent to the restriction of the standard Euclidean metric on
${\bf R}^{n + 1}$ to ${\bf S}^n$.

        Suppose now that $M$ is some kind of nice surface in ${\bf
R}^n$, which may be of any dimension.  Suppose also that $M$ is
connected in the sense that for every pair of points $x$, $y$ in $M$
there is a continuously-differentible path on $M$ that goes from $x$
to $y$.  In this situation, it is natural to try to define the
distance between $x$ and $y$ to be the length of the shortest curve on
$M$ that goes from $x$ to $y$.  In particular, it can be shown that
such a curve exists under suitable conditions.  Alternatively, one can
avoid the issue by defining the distance from $x$ to $y$ to be the
infimum or greatest lower bound of the lengths of the paths on $M$
that go from $x$ to $y$.  It is easy to see that this automatically
satisfies the triangle inequality, using the fact that a path from $x$
to $y$ may be combined with a path from $y$ to $z$ to get a path from
$x$ to $z$ whose length is equal to the sum of the lengths of the
other two paths for any $x, y, z \in M$.  If $M$ is the unit sphere,
then it is well known that this is the same as the metric described in
the previous paragraph.

        Consider the space $C([0, 1])$ of continuous real-valued
functions on the unit interval $[0, 1]$.  One way to define a metric
on $C([0, 1])$ is by
\begin{equation}
        d_1(f, g) = \int_0^1 |f(x) - g(x)| \, dx.
\end{equation}
Another way is to use
\begin{equation}
        d_\infty(f, g) = \max_{0 \le x \le 1} |f(x) - g(x)|.
\end{equation}
Note that the maximum of $|f(x) - g(x)|$ is attained on $[0, 1]$, by
the extreme value theorem.  It is not difficult to check that $d_1(f,
g)$ and $d_\infty(f, g)$ both satisfy the requirements of a metric on
$C([0, 1])$, and that
\begin{equation}
        d_1(f, g) \le d_\infty(f, g)
\end{equation}
for every $f, g \in C([0, 1])$.

\subsection*{The unit sphere}

        Let us return to the distance function $d_{{\bf S}^n}(x, y)$
defined on the unit sphere ${\bf S}^n$ in ${\bf R}^{n + 1}$ as before.
Note that $d_{{\bf S}^n}(x, y)$ attains its maximal value $\pi$ when
$x$ and $y$ are antipodal, which means that $x = -y$.  Let $\|x\| =
\Big(\sum_{j = 1}^{n + 1} x_j^2\Big)^{1/2}$ be the standard Euclidean
norm of $x = (x_1, \ldots, x_{n + 1}) \in {\bf R}^{n + 1}$, so that
$\|x - y\|$ is the standard Euclidean metric on ${\bf R}^{n + 1}$.
The spherical distance on ${\bf S}^n$ is related to the standard
Euclidean metric on ${\bf R}^{n + 1}$ by
\begin{equation}
\label{sin (frac{d_{{bf S}^n}(x, y)}{2}) = frac{||x - y||}{2}}
 \sin \Big(\frac{d_{{\bf S}^n}(x, y)}{2}\Big) = \frac{\|x - y\|}{2},
\end{equation}
which holds for every $x, y \in {\bf S}^n$.  This can be verified by
considering the line through $0$ and the mid-point $(x + y)/2$ of the
line segment that connects $x$ and $y$, which are perpindicular to
each other.

        If $x, y, z \in {\bf S}^n$, then we would like to show that
\begin{equation}
\label{d_{{bf S}^n}(x, z) le d_{{bf S}^n}(x, y) + d_{{bf S}^n}(y, z)}
        d_{{\bf S}^n}(x, z) \le d_{{\bf S}^n}(x, y) + d_{{\bf S}^n}(y, z).
\end{equation}
We may as well suppose that $x \ne y \ne z$ and $d_{{\bf S}^n}(x, y) +
d_{{\bf S}^n}(y, z) < \pi$, since this is trivial otherwise.  In
particular, $d_{{\bf S}^n}(x, y), d_{{\bf S}^n}(y, z) < \pi$, which
implies that $x, z \ne -y$.  Let $P(x, y)$ be the $2$-dimensional
plane in ${\bf R}^{n + 1}$ passing through $x$, $y$, and $0$, and let
$C(x, y)$ be the circle which is the intersection of $P(x, y)$ with
${\bf S}^n$, as before.  If $z$ is an element of $C(x, y)$, then
(\ref{d_{{bf S}^n}(x, z) le d_{{bf S}^n}(x, y) + d_{{bf S}^n}(y, z)})
is clear.

        Consider
\begin{eqnarray}
 \Sigma & = & \{w \in {\bf S}^n : \|y - w\| = \|y - z\|\}  \\
        & = & \{w \in {\bf S}^n : d_{{\bf S}^n}(y, w) = d_{{\bf S}^n}(y, z)\},
                                                          \nonumber
\end{eqnarray}
which is a sphere of dimension $n - 1$ in ${\bf R}^{n + 1}$.  The
intersection of $\Sigma$ with $C(x, y)$ consists of exactly two points
$u$, $v$, and we can label them in such a way that
\begin{equation}
        \|x - u\| \le \|x - v\|,
\end{equation}
which is equivalent to
\begin{equation}
        d_{{\bf S}^n}(x, u) \le d_{{\bf S}^n}(x, v).
\end{equation}
The main point now is that
\begin{equation}
\label{||x - u|| le ||x - w|| le ||x - v||}
        \|x - u\| \le \|x - w\| \le \|x - v\|
\end{equation}
for every $w \in \Sigma$, which implies that
\begin{equation}
        d_{{\bf S}^n}(x, u) \le d_{{\bf S}^n}(x, w) \le d_{{\bf S}^n}(x, v).
\end{equation}
In particular, we can apply this to $w = z$, to get that
\begin{equation}
        d_{{\bf S}^n}(x, z) \le d_{{\bf S}^n}(x, v).
\end{equation}
This implies (\ref{d_{{bf S}^n}(x, z) le d_{{bf S}^n}(x, y) + d_{{bf
S}^n}(y, z)}), since we already know that (\ref{d_{{bf S}^n}(x, z) le
d_{{bf S}^n}(x, y) + d_{{bf S}^n}(y, z)}) holds when $z \in C(x, y)$.
To get (\ref{||x - u|| le ||x - w|| le ||x - v||}), note that $\Sigma$
is centered at a point $\sigma$ on the line segment connecting $y$ to
$-y$, and that $\Sigma$ is contained in the $n$-dimensional plane $H$
in ${\bf R}^{n + 1}$ that passes through $\sigma$ and is perpindicular
to the line $L(y)$ through $y$ and $0$.  Let $L'$ be the line
contained in $P(x, y)$ that passes through $\sigma$ and is
perpindicular to $L(y)$, which is the same as the intersection of
$P(x, y)$ with $H$.  Also let $x'$ be the orthogonal projection of $x$
in $L'$, so that $x' \in L'$ and the line $L''$ passing through $x$
and $x'$ is perpindicular to $L'$.  If $x$ is already an element of
$L'$, then $x = x'$.  Note that $L''$ is parallel to $L(y)$, and that
$x'$ is the same as the orthogonal projection of $x$ in $H$.  Thus
\begin{equation}
        \|x - w\|^2 = \|x - x'\|^2 + \|x' - w\|^2
\end{equation}
for every $w \in H$, so that maximizing or minimizing $\|x - w\|$ for
$w \in \Sigma$ is the same as maximizing or minimizing $\|x' - w\|$ on
$\Sigma$.  The maximum and minimum of $\|x' - w\|$ for $w \in \Sigma$
are attained on the line $L'$, since $L'$ passes through $x'$ and the
center $\sigma$ of $\Sigma$ and is contained in $H$.  We also have that
\begin{equation}
        L' \cap \Sigma = P(x, y) \cap \Sigma = C(x, y) \cap \Sigma,
\end{equation}
where the first step uses the fact that $L' = P(x, y) \cap H$, the
second step uses the fact that $C(x, y) = P(x, y) \cap {\bf S}^n$, and
both steps use the fact that $\Sigma = {\bf S}^n \cap H$.  It follows
that the maximum and minimum of $\|x' - w\|$ on $\Sigma$ are attained
on $C(x, y) \cap \Sigma$, which are the same as the maximum and
minimum of $\|x - w\|$ on $\Sigma$, as desired.

\subsection*{Infinite series}

        If $x$ is a real number, then it is well known and easy to see that
\begin{equation}
        (1 - x) \sum_{j = 0}^n x^j = 1 - x^{n + 1}
\end{equation}
for each nonnegative integer $n$.  Here $x^n$ is interpreted as being
equal to $1$ when $n = 0$, even when $x = 0$.  In particular,
\begin{equation}
        \sum_{j = 0}^n x^j = \frac{1 - x^{n + 1}}{1 - x}
\end{equation}
when $x \ne 1$.  If $|x| < 1$, then $|x^n| = |x|^n \to 0$ as $n \to
\infty$, and so
\begin{equation}
        \lim_{n \to \infty} \sum_{j = 0}^n x^j = \frac{1}{1 - x}.
\end{equation}
Equivalently, $\sum_{j = 0}^\infty x^j$ converges as an infinite
series of real numbers, and
\begin{equation}
        \sum_{j = 0}^\infty x^j = \frac{1}{1 - x}.
\end{equation}

        Suppose now that $p$ is a prime number, $x$ is a rational
number, and $|x|_p < 1$.  This implies that $|x^n|_p = |x|_p^n \to 0$
as $n \to \infty$, which means that $x^n \to 0$ as $n \to \infty$ with
respect to the $p$-adic metric on ${\bf Q}$.  Similarly,
\begin{equation}
        \lim_{n \to \infty} \sum_{j = 0}^n x^j = \frac{1}{1 - x}
\end{equation}
with respect to the $p$-adic metric on ${\bf Q}$.  This is the same as
saying that $\sum_{j = 0}^\infty x^j$ converges as an infinite series
of rational numbers with respect to the $p$-adic metric, with the same
sum $1/(1 - x)$ as before.  For example, one can apply this to $x = p$,
to get that $\sum_{j = 0}^\infty p^j$ converges with respect to the
$p$-adic metric, with sum equal to $1/(1 - p)$.

        If an infinite series $\sum_{j = 1}^\infty a_j$ of real
numbers converges, then it is well known that $a_j \to 0$ as $j \to
\infty$ with respect to the standard metric on ${\bf R}$.  It is also
well known that $\sum_{j = 1}^\infty a_j$ may not converge even though
$\lim_{j \to \infty} a_j = 0$, e.g., when $a_j = 1/j$.  By contrast,
it can be shown that an infinite series $\sum_{j = 1}^\infty a_j$ of
$p$-adic numbers converges if and only if $a_j \to 0$ as $j \to
\infty$ with respect to the $p$-adic metric.  In both cases,
convergence of an infinite series is defined to mean convergence of
the corresponding sequence of partial sums, which is equivalent to
asking that the sequence of partial sums be a Cauchy sequence, by
completeness.  In the $p$-adic case, it is much easier to check that a
sequence is a Cauchy sequence, because of the stronger form of the
triangle inequality.

\subsection*{Norms on ${\bf R}^n$}

        A real-valued function $\|x\|$ on ${\bf R}^n$ is said to be a
\emph{norm} if it satisfies the following three conditions.  First,
$\|x\| \ge 0$ for every $x \in {\bf R}^n$, with $\|x\| = 0$ if and
only if $x = 0$.  Second,
\begin{equation}
\label{||t x|| = |t| ||x||}
        \|t \, x\| = |t| \, \|x\|
\end{equation}
for every $t \in {\bf R}$ and $x \in {\bf R}^n$, where
$t \, x = (t \, x_1, \ldots, t \, x_n)$.  Third,
\begin{equation}
\label{||x + y|| le ||x|| + ||y||}
        \|x + y\| \le \|x\| + \|y\|
\end{equation}
for every $x, y \in {\bf R}^n$, where $x + y = (x_1 + y_1, \ldots, x_n
+ y_n)$.  If $\|x\|$ is a norm on ${\bf R}^n$, then it is easy to see that
\begin{equation}
        d(x, y) = \|x - y\|
\end{equation}
defines a metric on ${\bf R}^n$.

        For example,
\begin{equation}
        \|x\|_2 = \Big(\sum_{j = 1}^n x_j^2\Big)^{1/2}
\end{equation}
is the standard Euclidean norm on ${\bf R}^n$, for which the
corresponding metric is the standard Euclidean metric $d_2(x, y)$,
defined earlier.  Similarly, it is easy to see that
\begin{equation}
        \|x\|_1 = \sum_{j = 1}^n |x_j|
\end{equation}
and
\begin{equation}
        \|x\|_\infty = \max_{1 \le j \le n} |x_j|
\end{equation}
are also norms on ${\bf R}^n$, and that these norms correspond to the
metrics $d_1(x, y)$ and $d_\infty(x, y)$ that were mentioned at the
beginning.  

        A subset $E$ of ${\bf R}^n$ is said to be \emph{convex} if for
every $x, y \in E$ and $t \in {\bf R}$ with $0 \le t \le 1$, we have that
\begin{equation}
        t \, x + (1 - t) \, y \in E.
\end{equation}
This is the same as saying that the line segment in ${\bf R}^n$
connecting $x$ to $y$ is contained in $E$ for every $x, y \in E$.
Let $\|x\|$ be a norm on ${\bf R}^n$, and let
\begin{equation}
        B = \{x \in {\bf R}^n : \|x\| \le 1\}
\end{equation}
be the closed unit ball in ${\bf R}^n$ corresponding to $\|x\|$.  It
is easy to see that $B$ is convex.  Moreover, $B$ is also
\emph{symmetric} about the origin in ${\bf R}^n$, in the sense that
$-x \in B$ for every $x \in B$.

        The definition of a norm makes sense on any vector space $V$
over the real numbers, such as the space $C([0, 1])$ of continuous
real-valued functions on the unit interval.  Any norm on $V$ determines
a metric on $V$ in the same way as before, and the closed unit ball
associated to the norm is convex and symmetric about the origin in $V$.

\subsection*{Distances on graphs}

        A \emph{graph} $G$ can be described by a set $V$ of
\emph{vertices}, and a subset $E$ of the set of subsets of $V$ with
exactly two elements.  Thus we say that there is an \emph{edge}
between two distinct elements $v$, $w$ in $V$ when $\{v, w \} \in E$.
A finite sequence $v_0, v_1, \ldots, v_n$ of vertices in $G$ is said
to determine a \emph{path} in $G$ if there is an edge between $v_j$
and $v_{j + 1}$ for each $j = 0, \ldots, n - 1$, in which case the
length of this path is defined to be $n$.  We say that $G$ is
\emph{connected} if for every pair of vertices $v, w \in V$ there is a
path $v_0, \ldots, v_n$ with $v_0 = v$ and $v_n = w$.  If $v = w$,
then we can simply take $n = 0$ and $v_0 = v = w$.

        If $G$ is a connected graph and $v, w \in V$, then the
distance $d(v, w)$ between $v$ and $w$ in $G$ may be defined to be the
smallest nonnegative integer $n$ for which there is a path in $G$ of
length $n$ that goes from $v$ to $w$.  It is easy to see that this
defines a metric on $V$.  In particular, if $v, w, z \in V$, then any
path in $G$ from $v$ to $w$ can be combined with a path from $w$ to
$z$ to get a path from $v$ to $z$, whose length is the sum of the
lengths of the paths from $v$ to $w$ and from $w$ to $z$.  This
implies that this definition of distance satisfies the triangle
inequality.

        This definition can be extended so that the distance between
adjacent vertices is any positive real number, depending on the vertices.
The edges can also be represented by segments or other curves, and included
in the associated metric space.  In this case, one can consider continuous
paths in the metric space, and not just discrete paths.

        Metric spaces are often discussed in books on basic analysis
and topology, in connection with the theory behind calculus.  A few
references along these lines are given below.

\end{document}